# MCMC Sampling of Directed Flag Complexes with Fixed Undirected Graphs

Florian Unger[*], Jonathan Krebs[†]

September 5, 2023

## Abstract

Constructing null models to test the significance of extracted information is a crucial step in data analysis. In this work, we provide a uniformly sampleable null model of directed graphs with the same (or similar) number of simplices in the flag complex, with the restriction of retaining the underlying undirected graph. We describe an MCMC-based algorithm to sample from this null model and statistically investigate the mixing behaviour. This is paired with a high-performance, Rust-based, publicly available implementation. The motivation comes from topological data analysis of connectomes in neuroscience. In particular, we answer the fundamental question: are the high Betti numbers observed in the investigated graphs evidence of an interesting topology, or are they merely a byproduct of the high numbers of simplices? Indeed, by applying our new tool on the connectome of C. Elegans and parts of the statistical reconstructions of the Blue Brain Project, we find that the Betti numbers observed are considerable statistical outliers with respect to this new null model. We thus, for the first time, statistically confirm that topological data analysis in microscale connectome research is extracting statistically meaningful information.

# 1 Introduction

## 1.1 Motivation

Inspecting spike correlation data, i.e. groups of neurons preferably spiking together, is a natural approach to investigate information processing in the brain. As the data comes in the form of a simplicial complex, the choice of using topological data analysis (TDA) for in-depth inspection is natural. This approach has seen recent success in [7, 18], but is severely limited by the data available: acquiring large-scale, yet detailed, in-vivo neuronal activity records remains infeasible, which restricts analysis to experiments of small scale or to in-silico simulations.

[*]Technische Universität Graz, Austria, florian.unger@igi.tugraz.at
[†]Friedrich-Alexander-Universität Erlangen-Nürnberg, Germany, jonathan.krebs@fau.de

**Why investigate connectomes?** Analyzing connectomes, i.e. the directed adjacency graph on the level of single neurons and synaptic connections between them instead, might mitigate these shortcomings: not only do they offer additional information due to the directionality added, it also makes massive amounts of data readily available. One kind would be recently developed large-scale sparse connectome reconstructions [2, 12]. These reconstructions stochastically generate connectomes from a large — but in comparison to biology still limited — set of biological facts (e.g. by using the position and orientation of different neuron types and their respective probabilities to form synaptic links), resulting in networks which match these biological statistics. Even more excitingly, there has been rapid progress in the field of dense reconstructions [9, 13, 17, 22]. In this approach the connectome of small volumes of actual brain tissue is determined. Electron microscopes together with advanced data processing methods achieve a resolution which allows one not only to map neurons (as with traditional light-based microscopy, which is unable to pick up synapses) but even allows to locate and count synaptic vesicles which are several magnitudes smaller. This scientific and technical progress enables an accurate and detailed reconstruction of networks of neurons and synapses far beyond researchers wildest dreams just a mere decade ago. With this kind of indisputable biological ground truth data becoming increasingly available, the accompanying development of analytical methods is most desirable. As it has seen success in the related analysis of spike correlation data, one of the potential methods trying to establish itself is TDA.

**Connectomic analysis with TDA.** When thinking about graphs through the lens of TDA, it is not immediately clear how to go about extracting complex information: in contrast to simplicial complexes gained by activity correlation, connectomic graphs themselves are just a simplicial complex of maximum dimension one (as vertices are of dimension 0 and edges are of dimension 1) and are therefore topologically not very interesting.

Thus one approach to enable extraction of information is to first construct the so-called flag complex. This represents a very specific case of a simplicial complex which contains *all* simplices compatible with the given graph. The space generated by this approach has a geometrical realisation and could be interpreted as the space describing potential signal information flow in the brain, where higher-dimensional simplices represent robust, dense, yet directed paths.

Once the flag complex has been computed by greedily searching for all simplices present, further topological analysis is possible, e.g. the calculation of Betti numbers in different dimensions in order to get an idea of the topology of the entire structure [15]; $q$-analysis, i.e. looking for connected simplicial pathways [16]; or the simplex completion probability [19], which endeavours to shed light on simplex generation mechanisms.

The above approach was pioneered in [15], where, amongst others, the dense reconstruction of the neural network of the roundworm Caenorhabditis Elegans (C. Elegans) [4, 20] was investigated (cited in Figure 1, left). Here, two findings stand out: compared to an Erdős–Rényi graph (ER-graph) of equal size and edge density, the connectome of C. Elegans seems to have, in all higher dimensions,

- more simplices *and*
- higher Betti numbers.

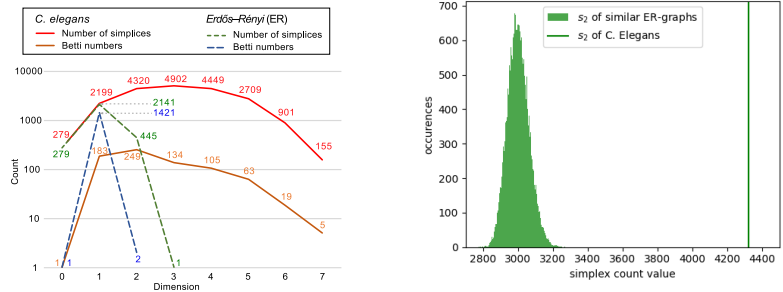

Figure 1: **Left:** Simplex Count and Betti numbers of C. Elegans compared to an ER-graph. Figure taken from [15]. **Right:** Histogram of the simplex count in dimension 2 of random ER-graphs with the same underlying undirected graph as C. Elegans.

But are these really two independent findings? Nontrivial homology expressed by nonzero or high Betti numbers in a given dimension requires, *by definition*, simplices in this dimension, which the ER-graph used for comparison does not provide. Undoubtedly aware of the problem, the authors of [15] never claimed to have found statistically interesting topology.

In spite of the traction gained by aboves TDA-based approach of connectomic analysis, the same fundamental question remains:

*What kind of Betti numbers are to be expected in such simplex-rich graphs?*

## 1.2 The Need for Sampleable Null Models

As a purely theoretical analysis currently seems to be intractable, we seek at least experimental, statistical insight. For that approach, one needs a suitable null model — in this case the set of all graphs with the identical simplex count seems appropriate. Whilst a parameterized null-model might be envisioned as well, we kept to a nonparametric version as in line with other research in network science. In particular, the condition of keeping the identical simplex counts enforces the same number of vertices and edges already, making this null model a more restricted version of the equally nonparametric ER-graphs so popular in network science.

Our desired insight into the expected topology may then be acquired straightforwardly by uniformly sampling a large number of graphs from that null model and analyzing their topology with the established toolchain.

**Sampling from simplex-rich graphs is difficult.** The problem lies in quick, efficient and uniform sampling of graphs from this set. Whilst constructing graphs with a specified number of simplices is straightforward, the resulting graphs are not uniformly distributed (or it would be very difficult to construct and prove).

Rejection sampling, i.e. the approach to first sample many ER-graphs and reject all those which don't have the right amount of simplices would work in principle, but is computationally infeasible. Connectomes are so simplex rich that they lie hundreds of standard deviations away from the expected mean and thus, in practice, every single graph would be rejected. This can be easily observed in Figure 1, right, where the distribution in a very similar scenario is depicted.

Another approach takes inspiration from stub matching algorithms, which maintain the degree sequence as required for the configuration model, yet another popular approach in network science. In this approach all edges are cut, leading to a stub per edge at each participating vertex. These stubs are then randomly paired, which results in retaining the number of edges *per vertex*, the so called degree sequence, of the original graph. Adapted for simplicial complexes (which are temporarily interpreted as hypergraphs), this approach has the potential to rapidly sample directed simplicial complexes with a comparable number of maximal simplices. Yet, for us it is unsuitable (see next paragraph), as the resulting complexes are very, very far from being flag[1]. Thus rejection sampling the results becomes, once again, computationally infeasible.

With no easy solution in sight, it is prudent to first investigate preceding research.

**Related Work.** Due to the novelty of the construction of directed flag complexes, there are no competing null models for this specific scenario. The closest established null model to the one we envision here is the so called Simplicial Configuration Model [23]. It is an extension of the configuration model, which maintains the degree sequence of an undirected graph (i.e. each vertex maintains the same number of edges). As such, it maintains the number of maximal simplices and additionally for each node the number of maximal simplices attached to it. While useful in many other scenarios (e.g. the initially mentioned [18]), this model does not satisfy our needs: besides being undirected and focusing on maximal simplices, it would *not* respect the property that our original simplicial complex is flag. The last difference especially is crucial, as it is already known that flag complexes are really specific edge cases of simplicial complexes. Thus any statistical difference between the connectomes and this null model might be attributed due the samples not being flag.

Nonetheless, their work demonstrated how a restrained random walk with small and potentially erroneous steps could be interpreted as a Markov Chain Monte Carlo approach which *guarantees*, under some reasonable assumptions, uniform sampling.

## 1.3 The Proposed Approach

**Restricting the null model further.** Straightforwardly adapting [23] turned out to be unviable: the difficulty lies not only with the additional directionality, but primarily with the additional complexity introduced with the upward-closure of flag-complexes (see Remark 2.7).

In order to start somewhere, we decided to conceptually factorize modifications of directed graphs into two distinct ways: connectivity and directionality. This allows us to treat the two resulting subproblems separately, resulting in the original problem becoming more approachable. It further allows for more detailed insight of the origin of the topological anomaly, should one be present. On the other hand, should one isolated aspect exhibit significant nonrandom behaviour already, one would expect that to still be visible in the general model.

---

[1] If one is primarily interested in simplicial complexes and does not care about the flag property, this approach is quite interesting, though. It is thus currently investigated in a similar project.

The connectivity is an undirected property. The notion of simplices reduces to the notion of cliques and the flag complex becomes the clique complex. Our previously envisioned null model would then require graphs with the same number of cliques. We believe that the approach used in [23] could be adapted to this model, although that has yet to be investigated in detail.

The directionality *without* the connectivity then respects a given undirected graph (i.e. the information whether $A$ and $B$ are connected at all), and is flexible only with respect to the direction of edges and the position of reciprocated edges (later referred to as double edges). In this work, the focus lies exclusively on directionality.

**Proposed Null Model** We thus propose a null model which enables uniform sampling of graphs that, with respect to a given graph, have:

- the same number of vertices and directed edges,
- a similar amount of simplices in the flag complex,
- the same underlying undirected graph.

Here 'similar' means an allowed, user-adjustable deviation from the simplex count, e.g. "within 1% of the original graph". The deviation is not meant as a parameter, more like a small error deemed acceptable, as maintaining the *precise* simplex count is extremely difficult.

**Overview of the utilized approach.** Even with the additional restriction of retaining the undirected graph, rejection sampling as explained above remains infeasible (see Figure 1, right). So far our research indicates that approaches on the level of simplicial complexes do not maintain the flag property. Thus, as the flag complex is completely determined by the graph, we have to modify the graph itself, then compute the flag complex and observe the resulting changes. For this approach, it is crucial that local changes in the graph result only in local changes in the flag complex. Exploiting this locality property is what keeps our approach computationally feasible.

Inspired by [23], we utilize, superficially, the same technique: a restricted random walk with resampling on the set of desired graphs. After a sufficient number of random steps, the point we end up at is uniformly distributed on the set of desired graphs. This is under the assumption that each point is as likely to enter as to exit (*double stochastic*) and all the desired points are connected (*irreducible*).

Our small steps conceptually consist of the so-called *simple moves*: flipping edges and moving double edges (depicted in Figure 3). The edges to be modified are chosen uniformly at random, and independent of previous steps, making the process a Markov Chain. As every change is reversible with the same probability, we automatically get double stochasticity. These moves do have the potential to change the simplex count though, which over time might accumulate, and subsequently lead to a loss of the desired property of a similar simplex count.

We thus reject changes which would lead the random walk too far astray and simply resample the current point instead (see Figure 2, red and blue arrows).

A common problem with restricted random walks is restricting too much, and thus splitting up the desired set into connected components, which are not interconnected by small steps anymore. In terms of Markov processes, this is expressed as the Markov Chain no longer being irreducible. To tackle this problem we utilize two tricks which extend the approach used in [23]. The first is rather general and could be used in similar settings: we allow the random walk to stray as far away from our originally desired set as necessary for it to become connected again (the yellow zone in Figure 2). This implies an additional, subsequent rejection sampling step, which is necessary to filter out undesired graphs, but incurs significant additional computational cost. Second, we allow for larger modifications, the so called complex moves (depicted in Figure 6 and 7). Though independently discovered, they bear quite some resemblance to the stub-shuffling approach mentioned above (although restricted to small parts of the graph). Besides helping with irreducibility, they have the additional benefit of speeding up the random walk by combining dozens of small, local steps into one big move. This not only speeds up computation, but also guarantees some stability of the simplex count as well.

The choice of a rather basic MCMC approach over more refined versions is based on the following observations: Hamilton-MCMC or the Wang-Landau algorithm and all derivations thereof are not applicable, as our setting is very discrete in nature. Metropolis-Hastings and Reversible Jump MCMC both coincide with our approach in this particular setting: the acceptance ratio of Metropolis-Hastings binarizes and thus boils down to our chosen strategy of resampling. Should the dimension of the space searched change at any point then reversible jumps are the method of choice (see e.g. [14]). But in the case investigated here the number of vertices and edges and thus the dimension is fixed, which negates any benefit of reversible jumps.

### 1.4 Structure of the Paper

In Section 2, we briefly recall basic definitions and results on graphs, homology, and MCMC-sampling. Furthermore, we describe general methods to adapt MCMC-sampling to sets only very indirectly describable, as is the case here.

In Section 3 we apply this general framework from the viewpoint of Markov Chains to our problem, first describing and analyzing very simple transition steps useful for a theoretical approach to the problem and then more complex ones, which have, amongst others, computational advantages. We further discuss boundary relaxations necessary to achieve irreducibility of the transition matrix.

In Section 4 we adopt the algorithmic standpoint and briefly discuss our high-performance, publicly available, Rust-based implementation. We then empirically investigate the mixing speed of this process by analyzing the Hamming distance over time and conducting a $\chi^2$ test for sample independence.

In Section 5 we apply this null-model to the connectome of C. Elegans, showing for the first time that this connectome seems to have an interesting, statistically abnormal topology compared to the null model designed and acquired in the previous sections. We further apply it to parts of the connectome of the Blue Brain Project, demonstrating that even large graphs with millions of edges are feasible with our implementation.

In Section 6 we discuss potential spin-offs, side problems, generalisations and the path to a null-model which is no longer connectivity-restrained.

## 2 Foundations

### 2.1 Graphs

The following definitions of graphs and subgraphs are primarily stated for self-containedness and to establish notation. Most readers may skip them if aware of the subject. Details of the definition of cliques in directed graphs (as *induced* complete subgraphs) and simplices (as not-necessarily induced complete, acyclic, oriented subgraphs) are necessary to follow Section 3 and Section 4 in all technical detail, though. The same holds for lemmata on the locality property of flag-complexes.

#### 2.1.1 Graphs, Cliques and Simplices

**Definition 2.1** (Directed and Undirected Simple Graphs)**.**

- A simple undirected graph is a pair $G = (V, E)$ of a finite set of vertices $V$ and a symmetric relation $E \subseteq (V \times V) \setminus \Delta_V$, where $\Delta_V = \{v, v \mid v \in V\}$. We refer to the set of simple undirected graphs by **Graph**. By excluding $\Delta_V$ from $E$ we exclude self-loops from $v$ to itself.

- A simple directed graph is a pair $G = (V, E)$ of a finite set of vertices $V$ and a relation $E \subseteq (V \times V) \setminus \Delta_V$, where $\Delta_V = \{(v, v) \mid v \in V\}$. We refer to the set of simple directed graphs by **DiGraph**. Here an edge $(i, j) \in E$ is interpreted as an edge from $i$ to $j$. We call an edge a single edge, when $(i, j) \in E$, but $(j, i) \notin E$. It is not forbidden that both $(i, j)$ and $(j, i)$ lie in $E$, in that case we call it a double edge.

- We define $\text{pr} : \textbf{DiGraph} \to \textbf{Graph}$, the projection from directed graphs to its underlying undirected graphs by

$$\text{pr}(V, E) = (\text{pr}(V), \text{pr}(E)) = (V, \{\{v, v'\} \mid (v, v') \in E\}).$$

- We call $\text{pr}(G)$ of a directed graph $G$ the underlying undirected graph.

- Let $S = (V_S, E_S)$ and $G = (V, E)$ be (directed or undirected) graphs with $V_S \subseteq V$ and $E_S \subseteq E$. Then $S$ is called a subgraph of $G$. Furthermore, if $E_S$ is maximal w.r.t $V_S$ and $E$, then $S$ is called the induced subgraph and is denoted with $G[S]$. In other words: a induced subgraph is only restricted in the set of vertices, but inherits all edges from its original graph.

A directed graph contains the additional information of edge direction. The projector pr forgets this directional information. It is surjective onto the set of all graphs, but not injective.

**Definition 2.2** (Cliques and Simplices)**.** Let $G = (V, E)$ be a directed graph.

- Let $V_c \subseteq V$ be a vertex set such that for all $i, j \in V_c$, $i \neq j$ there is $(i, j) \in E$ or $(j, i) \in E$. Then $c = G[V_c]$ is called a clique of order $n = |V_c|$, or an $n$-clique.

- Let $c$ be a clique which is not a subgraph of another clique. Then $c$ is called a maximal clique.

- A simplex in $G$ is a subgraph $\sigma = (V_\sigma, E_\sigma)$ of a clique together with a total order $\prec$ on $V_c$ represented by its edges, i.e. edges point from lower to higher vertices w.r.t. $\prec$: $E_\sigma = \{(i, j) \in E \mid i, j \in V_c \text{ and } i \prec j\}$. With $d + 1 = |V_c|$ we refer to $d$ as the dimension of $\sigma$ and call $\sigma$ a $d$-simplex. A simplex is often denoted by giving the vertices in rising order: $\sigma = [v_0, ... v_d]$ implies that that $v_i \prec v_j$ iff $i < j$.

Cliques are usually defined over undirected graphs, where their definition is analogous: $\text{pr}(c)$ is a clique in an undirected graph. Simplices could also be defined as complete, oriented, acyclic subgraphs of $G$, respectively oriented, acyclic subgraphs of a clique in $G$.

**Example 2.3.** A $d$-simplex is a subgraph of a $(d+1)$-clique. But a $(d+1)$-clique does not have to contain a $d$-simplex: indeed, there can be zero, one or even several $d$-simplices per $(d+1)$-clique. Here we observe three 3-cliques, all of them containing the 0-simplices $[0], [1]$ and $[2]$ and the 1-simplices given by the edges. The one to the left contains no 2-simplex, as it is cyclic. The one in the middle contains exactly one 2-simplex which is $[0, 1, 2]$. The one to the right even contains two 2-simplices: $[0, 1, 2]$ and $[0, 2, 1]$.

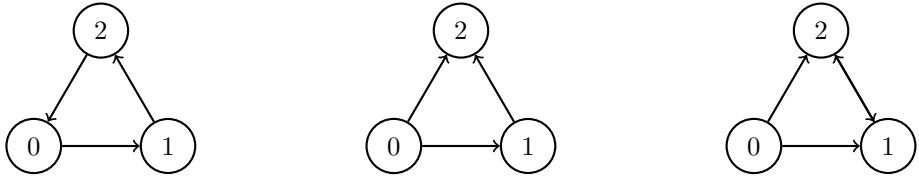

**Remark 2.4.** While Example 2.3 might, on first glance, imply that $d$-simplices in $(d+1)$-cliques are quite common, they are actually very rare in higher dimensions. Indeed, let $k = \frac{d^2+d}{2}$ be the number of edges in a $d$-simplex, leading us to $2^k$ potential (single-edge only)-configurations of the edges in a $(d+1)$-clique. At the same time there are "only" $(d + 1)!$ potential total orders on the set of vertices, leading to a ratio of $r(d) = \frac{(d+1)!}{2^{\frac{d^2+d}{2}}}$ of $d$-simplices per $(d+1)$-clique. That function shrinks quite quickly (values rounded):

| $d$    | 1    | 2   | 3     | 4     | 5    | 6     | 7      |
|--------|------|-----|-------|-------|------|-------|--------|
| $r(d)$ | 100% | 75% | 37.5% | 11.8% | 2.2% | 0.24% | 0.015% |

This effect is somewhat countered by occurrences of double edges and the rich number of lower-dimensional cliques in a high-dimensional clique. Still, there is a good chance that flipping a single edge destroys a high-dimensional simplex.

**Remark 2.5.** As Example 2.3 demonstrates, one clique may exhibit multiple simplices if the graph contains double edges. Indeed, an $n$-clique with $K = \binom{n}{2}$ double edges (i.e. every edge is a double edge) contains exactly $n!$ simplices of dimension *n*-1, as every order on the vertices could be expressed with the available edges. With *K*-1 double edges, the number of (*n*-1)-simplices is exactly $\frac{n!}{2}$. With two missing double edges the minimum amount of simplices drops to $\frac{n!}{4}$, whereas the maximum amount of simplices seems to follow the integer sequence A058298[2]. Once there are 3 double edges missing, one can enforce a 3-cycle, disabling all (*n*-1)-simplices within this *n*-clique. This leads to the following table:

| $k$ | 0 | 1 | 2 | 3 | 4 | 5 | 6 | ... | $K-3$ | $K-2$ | $K-1$ | $K$ |
|---|---|---|---|---|---|---|---|---|---|---|---|---|
| min | 0 | 0 | 0 | 0 | 0 | 0 | 0 | ... | 0 | $\frac{n!}{4}$ | $\frac{n!}{2}$ | $n!$ |
| max | 1 | 2 | 3 | 6 | 8 | 12 | 24 | ... | $\frac{n!}{4}$ | $\frac{n!}{3}$ | $\frac{n!}{2}$ | $n!$ |

### 2.1.2 Directed Flag Complexes and their Locality Property

**Definition 2.6** (Directed Flag Complexes and the Number of Simplices). Let $G = (V, E)$ be a simple directed graph and $D$ the dimension of the highest-dimensional simplex in $G$. The directed flag complex $d\mathcal{F}l(G)$ is a tuple of length $D+1$ containing in the respective position the set of all *d*-dimensional simplices of $G$.

We then define $s : \textbf{DiGraph} \to \mathbb{N}^{D+1}$, $s(G) = [s_0, s_1, ..., s_D]$ where $s_d$ yields the number of simplices of dimension $d$ in $d\mathcal{F}l(G)$.

**Remark 2.7.** Flag complexes are simplicial complexes, but the reverse is not necessarily true. A flag complex could be understood (or even defined) as the maximal simplicial complex over a given graph. This has a significant, and difficult to predict, impact on small modifications done to the complex: adding a high-dimensional simplex to a standard simplicial complex merely implies adding lower dimensional simplices as well, as every subset of a simplex must be contained (this property is referred to as downward closure). This is in contrast to flag complexes, where adding a (high-dimensional) simplex might not only require to add new lower-dimensional simplices, but then, as additional consequence of edges added due to the downward closure, new higher-dimensional simplices might appear in the graph and are added to the flag complex as well. One could refer to this property as upward-closure. This additional upward closure makes it difficult to adapt samplers designed for simplicial complexes to flag complexes.

The construction of flag-complexes is quite local: if there is a small local change in the directed graph, then the flag complex of that graph will remain the same for anything not in the direct vicinity of the changes.

**Definition 2.8** (Edge-Neighbourhoods). Let $G = (V, E)$ be a simple directed graph. Let $(i, j) \in E$. Then the neighbourhood of that edge is a subset of *V*:

$$\text{nbhd}((i,j)) = \{i, j\} \cup \Big( (\text{In}(i) \cup \text{Out}(i)) \cap (\text{In}(j) \cup \text{Out}(j)) \Big),$$

where $\text{In}(j) = \{i \in V \mid (i,j) \in E\}$ are the vertices with incoming edges and $\text{Out}(i) = \{j, \in V | (i,j) \in E\}$ are the vertices with outgoing edges.

This definition straightforwardly extends to multiple edges, with $D \subseteq E$:

$$\text{nbhd}(D) = \bigcup_{e \in D} \text{nbhd}(e).$$

---

[2]https://oeis.org/history?seq=A058298

**Lemma 2.9.** *Let $G = (V, E)$ and $G' = (V, E')$ be two directed graphs with the same underlying undirected graph. Let $D := (E' \setminus E) \cup (E \setminus E')$ be the (symmetric) difference between these two graphs. Then*

$$d\mathcal{F}l(G') = d\mathcal{F}l(G'[\text{nbhd}(D)]) \uplus (d\mathcal{F}l(G) \setminus d\mathcal{F}l(G[\text{nbhd}(D)]))$$

*Proof.* To keep the equations more readable, we shortcut $A := G[\text{nbhd}(D)]$ and $B := G'[\text{nbhd}(D)]$. We start with the following complement-defining equations:

$$d\mathcal{F}l(G) = d\mathcal{F}l(A) \uplus d\mathcal{F}l(A)^C$$
$$d\mathcal{F}l(G') = d\mathcal{F}l(B) \uplus d\mathcal{F}l(B)^C.$$

Our desired statement holds as soon as $d\mathcal{F}l(A)^C = d\mathcal{F}l(B)^C$ as we can then substitute $d\mathcal{F}l(B)^C$ with $d\mathcal{F}l(G) \setminus d\mathcal{F}l(A)$ in the second equation.

Let w.l.o.g. $\sigma \in d\mathcal{F}l(A)^C$. Then $\sigma$ is not a simplex in $d\mathcal{F}l(A)$ and thus $\text{pr}(\sigma)$ not a (undirected) clique in $\text{pr}(A)$. As the underlying undirected graphs of $A$ and $B$ are the same, $\text{pr}(\sigma)$ is also not a (undirected) clique in $\text{pr}(B)$. As a simplex is by definition also a clique in the underlying undirected graph, $\sigma$ cannot be a simplex in $d\mathcal{F}l(B)$.

Still, $\sigma \in d\mathcal{F}l(G')$ as all the edges of $\sigma$ are disjunct to $D$ and are thus both in $G$ and $G'$. Thus $\sigma \in d\mathcal{F}l(B)^C$. The same holds true for the other direction. □

This result might look technical, but is crucial for the implementation, as it allows for highly efficient calculations of updated simplex counts after local modifications of the graph. See Section 4.1.3 for the application.

## 2.2 Sampling via Markov Chain Monte Carlo

### 2.2.1 Unrestricted MCMC-Sampling

We recap the basics of Markov Chain Monte Carlo (MCMC) sampling, as it is the foundation of everything onward. For a concise introduction already in the context of graph sampling we refer to [5]. For a more exhaustive, general work with focus on mixing behaviour we refer to [8].

**Definition 2.10** (Time-Independent Markov Chains)**.** Let $T$ be a $N \times N$ stochastic matrix, i.e. it is non-negative and the sum of every row is 1. Let $\tau$ be the starting distribution over $S := \{1, \dots, N\}$. This forms a stochastic process, the time-independent Markov Chain $\mathcal{M}(S, T, \tau)$. It yields for time $t \in \mathbb{N}$ a random variables $X_t$ which is distributed according to $\tau T^t$.

One could describe $(S, T)$ as a weighted directed graph by interpreting $T$ as an adjacency matrix, but to avoid confusion with the graphs of our null model (which are points in $S$) we stick to the notion and notation of stochastic matrices.

**Definition 2.11** (Aperiodic, Doubly Stochastic, Irreducible)**.** Transition matrices may have the following properties:

**aperiodic** If there is a point $s \in S$ from which we start and return to at times $t_1, t_2, \dots$ then the greatest common divisor of $t_1, t_2, \dots$ shall be 1 for aperiodicity of $T$. On the level of graphs: the greatest common divisor of the length of all cycles shall be 1.

**doubly stochastic** The transpose of $T$ is stochastic as well, i.e. the sum over each row *and* each column is one. The corresponding graph property would be regularity.

**irreducible** For each $s, s' \in S$ there exists a $t \in \mathbb{N}$ such that $s'$ is reachable from $s$ in $t$ steps. In graph terms it would be strongly connected.

A matrix which is aperiodic and irreducible is called ergodic (as its Markov Chain is ergodic).

Ergodic Markov Chains have an equilibrium distribution, and if the transition matrix is doubly stochastic, this equilibrium distribution is uniform:

**Theorem 2.12** (Uniform Equilibrium Distribution)**.** *Let $S$ be a finite set of states and $T$ a transition matrix which is ergodic and doubly stochastic. Then the random variable $X_t$ defined by the Markov Chain $\mathcal{M}(S, T, \tau)$ converges, no matter the starting distribution $\tau$, to the uniform distribution, i.e. with $\nu = \frac{1}{|S|}$ and $\chi_t$ being the distribution of $X_t$:*

$$||\chi_t - \nu||_1 \to_{t \to \infty} 0.$$

This theorem has an important practical application: let us assume a random walk obeys the rules of an ergodic, doubly stochastic Markov Chain and starts somewhere in $S$. Then the theorem guarantees that the distribution it would end up with is approximately uniform. This approximation becomes more accurate the longer the random walk.

This process can be used to uniformly sample from $S$ (assuming an initial starting point is supplied). The approach is called *Markov Chain Monte Carlo Sampling*.

### 2.2.2 Mixing Time and Sampling Distance

Waiting till infinity is terribly boring. For practical applications concerned with uniformly sampling from $S$ the Markov Chain is usually run and subsampled every $k$ steps, where $k$ is called the sampling distance. Describing the (asymptotic) relationship between time and distance to the uniform distribution is referred to as studying the mixing behaviour, where rigorous theoretical results are often non-trivial. A more modest and applied approach — and thus the one pursued in this work — tries to estimate autocorrelation empirically.

### 2.2.3 Restricted MCMC-Sampling

Sometimes one strives to sample from a desired subset $A \subseteq S$. If $A$ can only be described in a very implicit way, the straightforward adaption of MCMC-sampling is to first try a step, with the same probabilities as before, and to check if one is still in the desired subset $A$. One has to be careful how to treat a misstep into $A^C$: simply rejecting a forbidden move would result in boundary effects on $A$ — boundary regions would be underrepresented as they are harder to reach (only from one side). The resulting distribution of the MCMC-sampling process would no longer be uniform, as the transition matrix would no longer be doubly stochastic.

A simple solution is the resampling-strategy: instead of merely rejecting a forbidden move and treating it as it has never happened, one resamples the current state instead.

This results in border regions which area hard to enter also being harder to leave, ensuring — at least up to the connected component — the desired uniform sampling quality.

**Definition 2.13** (Restricted Transition Matrix). Let $T$ be a transition matrix over a finite set of states $S$ and $A$ the subset of allowed states. Then with $a, a' \in A$ the restricted transition matrix $T_A$ is defined by

$$T_A(a, a') := \begin{cases} T(a, a') & \text{for } a \neq a' \\ T(a, a') + \sum_{s \in A^C} T(a, s) & \text{for } a = a'. \end{cases}$$

This trick is known in the context of degree-maintaining Markov Chains as "swap-and-hold" or "trial-swap" (see [1, 5]) and ensures that the resulting restricted transition matrix retains most of its MCMC-properties:

**Lemma 2.14.** *Let $T$ be a transition matrix which is aperiodic and doubly stochastic. Then the restricted transition matrix $T_A$ is aperiodic and doubly stochastic as well.*

Note that irreducibility is not so easily inherited — connected components of $A$, even if they were connected by $T$, are not necessarily connected by $T_A$ anymore (See Figure 2).

Indeed, this is the main concern throughout this paper. We developed the following general strategy to deal with it:

**Definition 2.15** ($A$-connecting subset). We call $R \subseteq S$ $A$-connecting iff $A$ is a subset of one irreducible component of $(R, T_R)$.

Instead of running the Markov Chain $\mathcal{M}(A, T_A, a)$ and potentially being stuck in one irreducible component of $A$, we run the Markov Chain on an $A$-connecting subset $R$, i.e. $\mathcal{M}(R, T_R, a)$ and subsequently rejection-sample the resulting samples for membership in $A$ (i.e. a sample $s$ is accepted if $s \in A$ and otherwise rejected).

As $A$ is a subset of the connected component of $R$, this results in uniformly sampled samples of $A$.

To optimize this approach one should strive for smaller, ideally even the minimum $A$-connecting subset. This is, in general, a hard problem and no universal solution is in sight.

# 3 Constructing the Transition Matrix

As motivated in the introduction, we pursue to sample uniformly from the following null model: let $G$ be a simple, directed graph with simplex count $s = s(G)$. We desire to sample from all graphs with the same underlying undirected graph, simplex counts close to $s$, i.e. between $s^- \leq s \leq s^+$ and exactly the same amount of vertices ($s_0$) and edges ($s_1$). With the formal definition introduced below, we refer to this set as $\mathcal{G}_{s^-}^{s^+}(G)$.

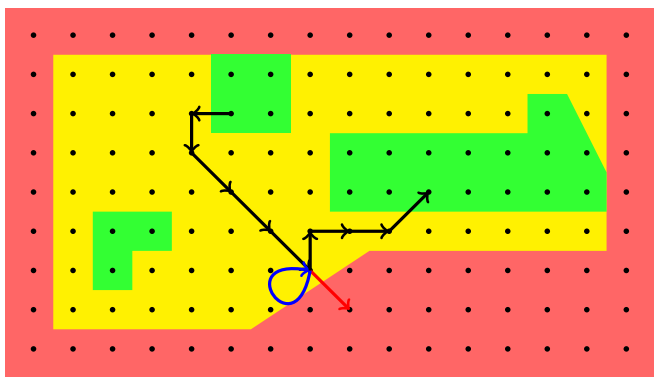

Figure 2: A sketch of a restricted random walk. Dots on the grid represent states in $S$. The green area represent the desired states in $A$, the yellow area represents a potential $A$-connecting subset $R$. The red area represents $R^C$ and is forbidden: Instead of entering the red area (red arrow) we immediately resample (blue arrow).

But how do we sample uniformly? In Section 2.2 we have introduced the theoretical background, restricted MCMC-processes with resampling and extended boundaries.

In this section we will discuss its application to our problem. This is primarily done by designing the required transition matrix. Here we show that it sticks to the desired subset $\mathcal{G}_{s^-}^{s^+}(G)$ and exhibits the desired properties of aperiodicity, double stochasticity and irreducibility necessary for an MCMC-process. Admittedly the last property is the most tricky one, instead of rigorous proofs we present a detailed discussion.

This is achieved by first investigating the so-called simple moves, a simple-to-understand and analyzable set of atomic moves which are enough to modify the directionality of a graph. They are complemented by the so-called complex moves, which have both theoretical and computational benefits by making use of the combinatorial nature of simplices. Both do not respect the simplex count by themselves though, so the random walk must be restricted as explained in the previous section. As these restriction might split $\mathcal{G}_{s^-}^{s^+}(G)$ into multiple disconnected subsets, the last part of this chapter discusses potential choices of connecting boundaries.

This section is complemented by Section 4, where everything is viewed from a more practical, algorithmic perspective and the mixing time is determined experimentally.

**Definition 3.1** (Simplex Count Boundaries). Let $G = (V, E)$ be an simple, directed graph and $D = \max\{|c| \mid c \text{ clique in } G\} - 1$ the maximal possible dimension of a simplex in $G$. Then $s^-, s^+ \in \mathbb{N}^D$ with $s_0^- = |V| = s_0^+$ and $s_1^- = |E| = s_1^+$ is a lower/upper simplex count boundary.

We denote with $\mathcal{G}_{s^-}^{s^+}(G)$ the set of all directed graphs $G'$ such that $\text{pr}(G') = \text{pr}(G)$ and $s^- \leq s(G') \leq s^+$ (for all dimensions simultaneously).

With $\mathcal{G}_0^\infty(G)$ we refer to the (not-really) boundaries that are 0 respectively $\infty$ for all dimensions higher than 1. This describes the set of all directed graphs with the same underlying graph as $G$ and the same number of (directed) edges, but no further restrictions on the simplex count.

## 3.1 Simple Moves

**Definition 3.2** (Simple Moves)**.** Let $G = (V, E)$ be a directed graph with at least one single edge (i.e. $(i,j) \in E$, but $(j,i) \notin E$). Then with $(i,j)$ a single edge we define the Single Edge Flip

$$\text{SEF}_{i,j} : G \mapsto G' \text{ where } G' = (V, E') \text{ with } E' = E \setminus \{(i,j)\} \uplus \{(j,i)\}.$$

Let further $G$ contain at least one double edge (i.e. both $(i,j)$ and $(j,i)$ are in $E$). With $(i,j)$ a single edge and $(k,l)$ a double edge we define the Double Edge Move

$$\text{DEM}_{i,j}^{k,l} : G \mapsto G' \text{ where } G' = (V, E') \text{ with } E' = E \setminus \{(k,l)\} \uplus \{(j,i)\}.$$

See Figure 3 for a graphical depiction.

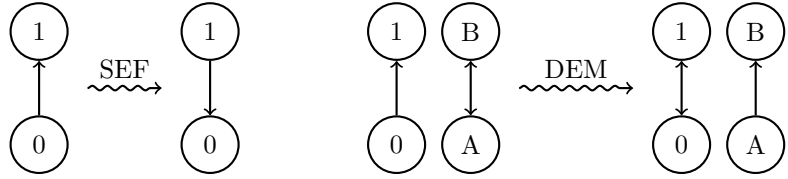

Figure 3: The simple moves single-edge-flip (SEF) and double-edge-move (DEM).

The following properties are straightforward to check and ensure that the image and coimage of SEF and DEM coincide, namely $\mathcal{G}_0^\infty(G)$:

**Lemma 3.3.** *Let $G$ be transformed to $G'$ by either SEF or DEM:*

1. *The underlying undirected graph is not modified:* $\text{pr}(G) = \text{pr}(G')$.

2. *The move retains the number of directed edges:* $s_1(G) = s_1(G')$.

3. $G' \in \mathcal{G}_0^\infty(G)$.

**Definition 3.4** (Unrestricted Simple Move Transition Matrix)**.** Let $G = (V, E)$ be a directed graph. We define the unrestricted transition matrices $T_{\text{SEF}}^{\text{u}}$ and $T_{\text{DEM}}^{\text{u}}$:

- $T_{\text{SEF}}^{\text{u}}(G, G')$ is given by the probability of sampling a single edge $(i,j) \in E$ by a uniform distribution over all single edges such that $G' = \text{SEF}_{i,j}(G)$. It is 0 if there is no such edge.

- $T_{\text{DEM}}^{\text{u}}(G, G')$ given by the probability of sampling a single edge $(i,j) \in E$ and a double edge $(k,l)$ uniformly (i.e. by a uniform distribution over the single respectively double edges) such that $G' = \text{DEM}_{i,j}^{k,l}(G)$. It is 0 if there are no such edges.

Furthermore, with $p_{\text{SEF}}, p_{\text{DEM}} \in (0,1)$ and $p_{\text{SEF}} + p_{\text{DEM}} = 1$, the unrestricted transition matrix of simple moves is the convex combination of each moves transition matrix:

$$T_{\text{SM}}^{\text{u}} = p_{\text{SEF}} T_{\text{SEF}}^{\text{u}} + p_{\text{DEM}} T_{\text{DEM}}^{\text{u}}.$$

As $V$ is finite, these matrices are finite.

**Lemma 3.5.** *Let $G$ be a simple, directed graph. The unrestricted transition matrix of simple moves on $G$, $T^u_{SM}$, is doubly stochastic and irreducible on $\mathcal{G}_0^\infty(G)$. If $G$ further contains at least one double edge, $T^u_{SM}$ is aperiodic as well.*

*Proof.* The previous lemma ensured that we do not leave $\mathcal{G}_0^\infty(U)$ with our moves. We check the remaining properties:

- Aperiodicity: The aperiodicity is due to the fact that there are cycles of length 2 and length 3 (see Figure 4), whose greatest common divisor is 1.

- Double Stochasticity: The move $\text{SEF}_{i,j}$ has the inverse $\text{SEF}_{j,i}$ and analogously, $\text{DEM}^{k,l}_{i,j}$ has the inverse $\text{DEM}^{j,i}_{k,l}$. The probability of picking the inverse move is the same as picking the original one: the number of single respectively double edges is retained through all moves and thus is the probability of picking one, as they're selected via a uniform distribution.

  This makes $T^u_{\text{SEF}}$ and $T^u_{\text{DEM}}$ symmetric and thus doubly stochastic. The same holds true for their convex combination $T^u_{\text{SM}}$.

- Irreducibility: We have to check that every graph $G' \in \mathcal{G}_0^\infty(G)$ is connected to $G$ by our set of moves. This can be done by defining a standard form which is reachable by all graphs in $\mathcal{G}_0^\infty(G)$: let us assume some arbitrary total order on the vertex set $V$. Then SEF and DEM could be used to flip all edges such that their ingoing vertex is ordered higher than their outgoing vertex (i.e. flip them to the upper right triangle of the adjacency matrix where possible) and push all double edges such that the tuple of outgoing vertex, ingoing vertex is minimal (i.e. push them as far top left as possible in the adjacency graph). This is a potential standard form reachable by all graphs in $\mathcal{G}_0^\infty(G)$ using only SEF and DEM. Thus $T^u_{\text{SM}}$ is irreducible.

□

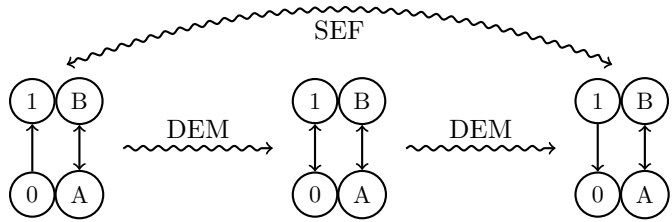

Figure 4: There is a cycle of length 3 (DEM from left to middle, DEM to the right, then SEF from right to left) and a cycle of length 2 (SEF from left to right and then immediately back).

## 3.2 Complex Moves

As we have witnessed e.g. in Figure 5 the simple moves SEF and DEM easily create and/or destroy simplices. As the formation of (especially) high-dimensional simplices by random chance is quite unlikely (Remark 2.4) we present here more complex moves to modify graphs. These combine many SEF and DEM at once while trying to retain $s(G)$ at the same time.

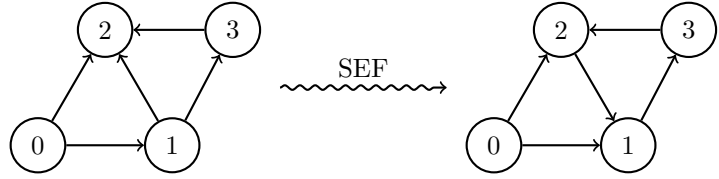

Figure 5: Moves potentially alter the simplex count: while the 2-simplex $[0, 1, 2]$ became the simplex $[0, 2, 1]$ by the flip of the edge $(1, 2)$, the former 2-simplex $[1, 3, 2]$ is lost.

This is achieved by applying permutations of the nodes onto the edges, i.e. $(i, j) \mapsto (\pi(i), \pi(j))$. That way the number of total orders on the vertices in a clique (which correspond to simplices) is retained. However, operations of that character could change the underlying undirected graph if applied to anything else than cliques.

**Definition 3.6** (Clique Permute). Let $G = (V, E)$ be a directed graph and $V_m$ a set of nodes such that $G[V_m] = (V_m, E_m)$ forms a maximal clique. Let $\pi \in S_{|V_m|}$ be a permutation on these vertices. Now

$$\mathrm{CP}_m^\pi : G \to G' = (V, E') \text{ where } E' = (E \setminus E_m) \uplus \pi(E_m),$$

with $\pi((i, j)) = (\pi(i), P(j))$ for $(i, j) \in E_m$.

See Figure 6 for a graphical depiction.

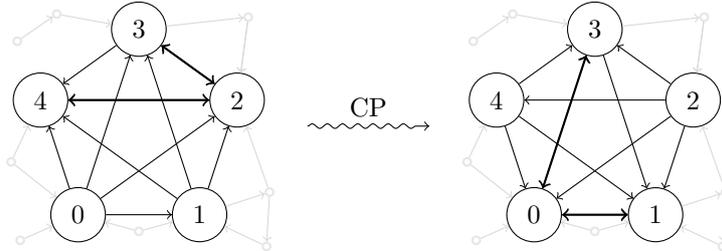

Figure 6: A demonstration of Clique Permute: The graph on the left exhibits the 4-simplices $[0, 1, 2, 3, 4]$ and $[0, 1, 3, 2, 4]$ over the same set of vertices. The application of Clique Permute with the permututation $\pi = (2, 4, 0, 3, 1)$ leads to the corresponding cliques $[1, 3, 2, 4, 0]$ and $[1, 3, 4, 2, 0]$. Small gray points and lines represent the graph outside of this clique, their respective connections are not modified.

Whereas Clique Permute applies within one clique, Clique Swap goes one step further: two equal-sized cliques swap all their edges. This not only permutes edges in two cliques, but has the potential to transfer many double edges at once, keeping the overall number of simplices highly depended on the number of double edges comparatively stable.

**Definition 3.7** (Clique Swap). Let $G = (V, E)$ be a directed graph. Let $V_1, V_2 \subseteq V$ be equal-sized sets of vertices such that they form, together with their induced edge sets $E_1$ and $E_2$, maximal cliques in $G$. Let further $\pi : V_1 \cup V_2 \to V_1 \cup V_2$ be a permutation with the additional property that $\pi(V_1) = \pi(V_2)$ (and thus also $\pi(V_2) = V_1$). With that we define the Clique Swap

$$\mathrm{CS}^\pi_{V_1, V_2} : G \mapsto G' = (V, E') \text{ where } E' = (E \setminus (E_1 \cup E_2)) \uplus \pi(E_1 \cup E_2),$$

where we define $\pi((i, j)) = (\pi(i), \pi(j))$ for $i, j \in V_1$ or $i, j \in V_2$.

The intended purpose and easiest case of the Clique Swap is to swap far-away cliques as sketched in Figure 7. However, these two cliques might overlap, requiring the additional technicalities above. On the other hand, if $V_1$ and $V_2$ overlap completely, we've just got the Clique Permute.

The requirement of the cliques to be maximal is not strictly necessary. In fact, by dropping it, CS is flexible enough to define SEF, DEM and CP at the same time. In our implementation, however, we stick to maximal cliques for a variety of reasons (see Section 4) and thus define them like above.

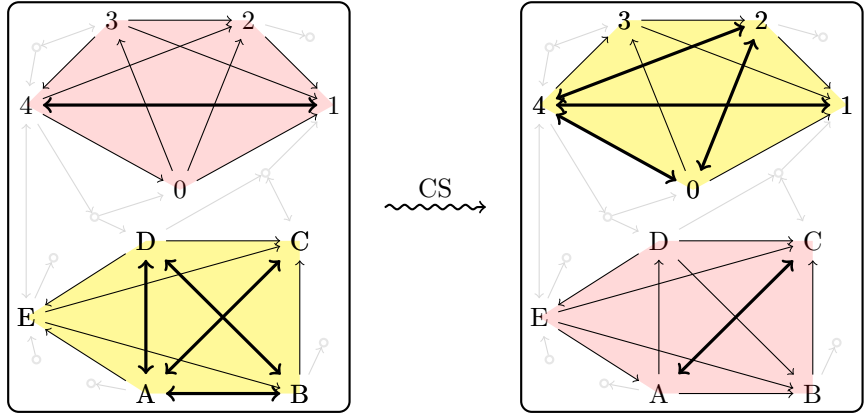

Figure 7: Clique Swap: On the left hand side we observe two 5-cliques with differing number of double edges. Clique Swap now swaps (via the bijection $(0, 1, 2, 3, 4) \mapsto (D,C,B,E,A)$) all the edges from the upper clique to the lower clique and vice versa. The remaining graph and their edges to these cliques remain untouched.

**Lemma 3.8.** Let $G = (V, E), V_1, V_2, \pi$ as in Definition 3.7 and $G' = \mathrm{CS}^\pi_{V_1, V_2}(G)$. Then

1. CS retains the number of directed edges: $|E| = |E'|$.

2. CS retains the underlying undirected graph: $\mathrm{pr}(G) = \mathrm{pr}(G')$.

3. $G' \in \mathcal{G}_0^\infty(G)$.

4. CS retains the number of simplices within $G[V_1] \cup G[V_2]$.

5. All of the above statements hold also true for $G' = CP^\pi_m(G)$.

*Proof.* We show for each

1. This is proven by showing that $|E_1 \cup E_2| = |\pi(E_1 \cup E_2)|$ and that $(E \setminus (E_1 \cup E_2)) \cap \pi(E_1 \cup E_2) = \emptyset$. The first straightforwardly follows from $\pi$ being a bijection: if $f$ is a bijection on $A$, then $f \times f$ is also a bijection on $A \times A$.

   For the second we note that $E \setminus (E_1 \cup E_2)$ contains no edges within $V_1$ and $V_2$ (it still potentially contains edges between $V_1$ and $V_2$, though!). $E_1 \cup E_2$ however contain *only* edges within $V_1$ respectively $V_2$. This property is retained by $\pi$ thanks to the requirement $\pi(V_1) = V_2$.

2. Consider w.l.o.g. $i, j \in V_1$. As $V_1$ forms a clique, $\{i, j\} \in \mathrm{pr}(E_1)$. But $\{i, j\} = \mathrm{pr}(\pi(\pi^{-1}(i), \pi^{-1}(j)))$, and thus $\pi^{-1}(i), \pi^{-1}(j) \in V_2$, which also forms a clique, implying that $\{\pi^{-1}(i), \pi^{-1}(j)\}$ must have been in $\mathrm{pr}(E_2)$. Thus $\mathrm{pr}(E_1 \cup E_2) = \mathrm{pr}(\pi(E_1 \cup E_2))$ and the desired property holds, as outside of $E_1 \cup E_2$ nothing was changed.

3. This is clear by the two previous statements.

4. A simplex is a clique together with directed edges corresponding to a local total order. As the undirected graph is retained, we must check only for the number of orders present. Let $\sigma = (v_0, v_1, \ldots, v_l)$ an $l$-simplex in $G[V_1] \cup G[V_2]$ and $\pi$ be a permutation with the required properties. Then the simplex $\sigma$ is transformed into the simplex $\pi(\sigma) = (\pi(v_0), \pi(v_1), \ldots, \pi(v_l))$, showing that no simplex is lost by CS. A similar argument with $\pi^{-1}$ shows that no simplex is gained as well. Thus the number of simplices is retained.

5. $\mathrm{CP}_m^\pi$ is actually a special case of CS where $V_1 = V_2$. Above results then follow immediately.

□

Note that CS does not retain the number of simplices in $G[V_1 \cup V_2]$.

We define the unrestricted transition matrix which comprises all the moves introduced so far:

**Definition 3.9** (Unrestricted Transition Matrix)**.** Let $G = (V, E)$ be a directed graph. Let $p_D$ be a distribution over all clique sizes occurring in $G$. We define the unrestricted transition matrices $T_{\mathrm{CP}}^{\mathrm{u}}$ and $T_{\mathrm{CS}}^{\mathrm{u}}$ the following way:

- $T_{CP}^{\mathrm{u}}(G, G')$ is defined by the product of the probabilities of picking a size $d$ (by $p_d$), then uniformly sampling a maximal clique $m$ and a permutation $\pi \in S_d$ such that $G' = \mathrm{CP}_m^\pi(G)$.

- $T_{CS}^{\mathrm{u}}(G, G')$ is defined by the product of the probabilities of picking a size $d$ (by $p_d$), then uniformly sampling two maximal cliques $V_1, V_2$ and a bijection as described in Definition 3.7 such that $G' = \mathrm{CS}_{V_1, V_2}^\pi(G)$.

Furthermore, with $(p_{\mathrm{SEF}}, p_{\mathrm{DEM}}, p_{\mathrm{CP}}, p_{\mathrm{CS}})$ forming a positive probability vector, the unrestricted transition matrix $T^{\mathrm{u}}$ is the convex combination of each moves transition matrix:
$$T^{\mathrm{u}} := \sum_{m \in \{\mathrm{SEF}, \mathrm{DEM}, \mathrm{CP}, \mathrm{CS}\}} p_m T_m^{\mathrm{u}}.$$

**Lemma 3.10.** *$T^u$ is symmetric and double stochastic.*

*Proof.* Maximal cliques are an undirected property of the graph, and as the underlying graph does not change, so neither $p_{M_d}$. Thus sampling $m_1$ and $m_2$ has the same chance at every time step (assuming $p_{M_d}$ does not change over time). As we picked our permutation $p_{\text{AB}}$ and $p_{\text{BA}}$ by a uniform distribution, picking the inverse permutation is just as likely. Thus the inverse move $G' \to G$ has the same probability as $G \to G'$ and the transition matrix $T^u_{\text{CS}}$ is symmetric. By definition $T^u_{\text{CS}}$ is stochastic. The same holds true for $T^u_{\text{CP}}$, and as a convex combination of symmetric stochastic matrices also for $T^u$. A symmetric stochastic matrix is double stochastic. □

**Corollary 3.11.** *Let $G = (V, E)$ be a directed graph with at least one double edge. Then $\mathcal{M}(\mathcal{G}_0^\infty(G), T^u, G)$ forms an ergodic Markov Chain with uniform equilibrium distribution.*

## 3.3 Relaxing Boundaries

The simple moves introduced do not strictly retain the number of simplices (see e.g. Figure 5). Even the more complex moves only retain the number of simplices *within* their cliques and not in their respective neighbourhood. To avoid straying to far from our desired simplex we must restrict the transition matrix $T^u$ as described in Section 2.2.3. This leads to an uncertainty if $\mathcal{G}_{s^-}^{s^+}$ is still connected. In this section we look for relaxations necessary to connect $\mathcal{G}_{s^-}^{s^+}$. We define the simplex-count equivalent of Definition 2.15:

**Definition 3.12** ($\mathcal{G}_{s^-}^{s^+}$-connecting relaxed bounds). We define the simplex count bounds $s^{--}$ and $s^{++}$ as $\mathcal{G}_{s^-}^{s^+}$-connecting with respect to an unrestricted transition matrix $T^u$, if for all $G, G' \in \mathcal{G}_{s^-}^{s^+}$ there exists a $t \in \mathbb{N}$ such that $(T_{s^{--}}^{s^{++}})^t_{G,G'} > 0$.

**Corollary 3.13.** *Let $G = (V, E)$ be a directed graph with at least one double edge. Let $s^-$, $s^+$ be simplex bounds and let $s^{--}$ and $s^{++}$ be $\mathcal{G}_{s^-}^{s^+}$-connecting bounds. Then $\mathcal{M}(\mathcal{G}_{s^{--}}^{s^{++}}(G), T^u, G)$ forms an ergodic Markov Chain with uniform equilibrium distribution.*

As we need to subsequently filter the output for membership in $\mathcal{G}_{s^-}^{s^+}$, we are interested in finding tight boundaries, meaning that the complement of $\mathcal{G}_{s^-}^{s^+}$ with respect to the irreducible component of $\mathcal{G}_{s^{--}}^{s^{++}}$ that contains $\mathcal{G}_{s^-}^{s^+}$ should not be unnecessarily large.

### 3.3.1 Single-Edge-Only Case

We first look at a special case: graphs that do not have any double edges. This is a property retained throughout all graphs with the same underlying undirected graph and the same number of edges and thus the moves introduced respect this property as well.

In this case we propose the following bounds:

**Conjecture 3.14.** *Let $G$ be a directed graph without double edges and $s^-, s^+$ simplex count boundaries. Then*

$$s^{--} := s^- \text{ and } s^{++} := \infty$$

*are $\mathcal{G}_{s^-}^{s^+}$-connecting.*

Assuming that we are looking for graphs with a simplex count far above the average, this is an easy and tight relaxed boundary. To convince ourselves of the tightness, we remind ourselves of the distribution of simplex counts as in Figure 1 (right) and Remark 2.4: While $s^{++} = \infty$ is probably not optimal (could be lower), it does not really matter, as with rising simplex count the number of graphs quickly decay, making it unlikely to sample them with a uniform sampler.

We now outlay our attempts to prove Conjecture 3.14: For easier analysis, we only take SEF into account.

The core idea is to connect $G, G'$ through star points which are directed acyclic graphs (DAGs). Let us assume there are DAGs $A$ and $A'$ connected to $G$ (respectively to $G'$). Now $A$ and $A'$ do not have any cycles and thus each correspond to a (global!) order on $V$. Thus both are connected by a series of SEF, as each SEF corresponds to a neighbouring transposition which are known to be generators of the symmetric group of permutations.

This leaves the question if and how $G$ is connected to its corresponding DAG $A$. This is equivalent to the question if *any* $G$ can be transformed to a DAG by a series of SEF, where each individual step may not decrease $s_2$ (every obstruction to a simplex in a clique is perceivable as a 3-cycle so eradicating them is enough).

In the case of $G$ being fully connected (i.e. a clique), we give, with Theorem A.3, a proof to this statement. The general case of non-complete graphs turns out to be surprisingly intricate to prove[3]: despite significant effort, we have to state this as the very elementary Conjecture A.4 for now. Yet, both Conjecture 3.14 and A.4 seem to be correct, as a computational search for counterexamples was unsuccessful so far, despite having checked billions of graphs.

Both Theorem A.3 and Conjecture A.4 also make a statement about the number of SEFs required: with them, we can bound the graph diameter of $(\mathcal{G}_{s^{--}}^{s^{++}}, T)$ to

$$|E| \leq \mathrm{diam}((\mathcal{G}_{s^{--}}^{s^{++}}, T)) \leq 2|E|,$$

a potentially useful statement for estimating mixing times.

### 3.3.2 Bounds for Graphs with Double Edges

In the previous subsection we have discussed suitable relaxations in the case of single-edge-only graphs. There the number of $d$-cliques was a natural and rather strict upper bound to connect everything.

This analysis is not enough to tackle the graphs motivating this work, biological connectomes. Indeed, C. Elegans contains 233 ( 11%) double edges and the Layers 1-4 of the BBP reconstruction contain 15258 ( 1%) double edges. Simultaneously, only a few double edges may result in a combinatorial explosion of the number of high dimensional directed simplices: an 8-clique in C. Elegans requires only $\binom{8}{2} = 28$ double edges to generate $8! = 40320$ 7-simplices. This questions the choice of the open upper boundary used in the single-edge-case: while it is still $(s^-, s^+)$-connecting, it is no longer a strict and thus efficient choice. An experiment confirmed that sampling from C. Elegans with an open upper boundary is not practical.

---

[3] In the words of one attendee of a week long open problem session where Conjecture A.4 was discussed: We admire the problem!

Thus we approach the double-edge-case differently. A rigorous treatment similar to the single-edge-only-case is currently out of scope, so we describe assumptions used in the implementation. They are based on the observation that the graphs investigated are quite large, yet somewhat regular and random.

- The graphs are large enough to have a considerable buffering-capacity: small changes in the simplex count caused by modifications in some parts of the graph can be (approximately) compensated in other, distant parts of the graph. It is then enough to ensure that boundaries are relaxed enough to allow for *a single* DEM as the graph may use its buffering capacity afterwards to return to the former number of overall simplices. This buffering capacity enables more complex local changes which would otherwise result in a (temporary) boundary violation if performed in direct succession.

- We can treat the simplex counts separately for each dimension: there are, for every dimension, enough maximal cliques available such that the lower dimensional effects of high dimensional changes can be compensated using solely the maximal cliques of these lower dimensions. But as these cliques are maximal the change will not propagate up again.

- In particular, there are many maximal 1-cliques which allow removing or adding double edges from the remaining graph by "hiding" them in these 1-cliques. This allows us to (locally) view DEM as removing or adding an edge.

With these assumptions we have to look for the worst-case scenario of connecting two graphs with similar simplex counts. Our approximation of a worst-case scenario is sketched in Figure 8: two *D*-sized cliques (*D* being maximal clique size), one (upper left) full with double edges and the other (lower left) void of any double edges. The upper clique thus exhibits $D!$ many simplices and the lower at most one. In the transformed graph (right) the double edges have switched their position. This is exactly the scenario CS would resolve in a single transition and was the original motivation to include this transition in the first place. CS entails modifying *multiple* edges simultaneously though, a violating critical to the buffering effect. Thus, for the theoretical connectivity analysis, we only consider DEM for the rest of this chapter.

Here we identify the following high-impact DEM. Removing the first edge of a clique full of double edges halves the number simplices in this clique (and analogously, adding the last one missing doubles the number of simplices). This effect is bounded primarily by the size of the clique, leading to $\Delta_d^{(1)} = \frac{1}{2}(d+1)!$ as a first approximation of a suitable boundary.

Usually, especially in the higher dimensions, $(d+1)!$ is (much) bigger than $s_d^+$, meaning that cliques full of double edges do not appear in the graphs desired. In these cases $\Delta_d^{(1)}$ is then an impractically generous boundary. A tighter one requires the integer sequence A058298 which describes in the $k$th entry the maximal number of $d$-simplices in a large enough clique with $k < (d+1)!$ double edges (see Remark 2.5). In the adapted situation of Figure 8 the maximal unavoidable simplex count loss or gain in a series of DEM is then (generously) bounded by
$$\Delta_d^{(2)} = \max_{k=1,\ldots,k'} \{\text{A058298}_k - \text{A058298}_{k-1}\}.$$

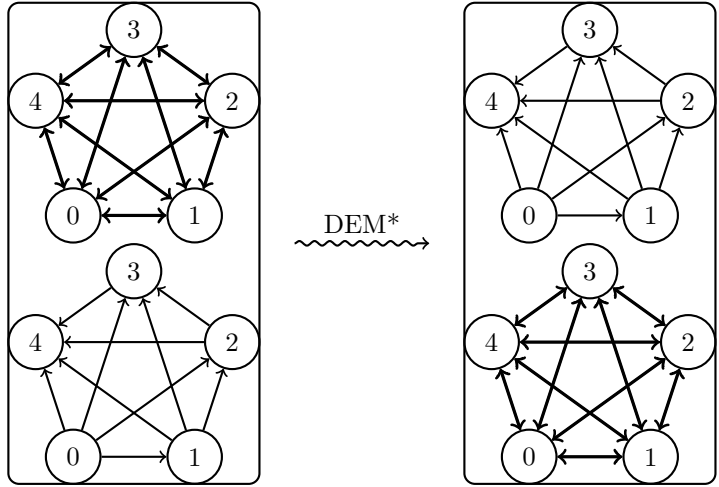

Figure 8: A potential worst-case scenario to find relaxed simplex bounds for.

It is a rather generous estimate as there are many cases where smaller changes in the clique are possible, e.g. by internally moving double edges instead of moving them out of the clique.

This handles the scenario depicted in Figure 8 and similar ones. This scenario is a bit optimistic though, as one edge may participate in many cliques and thus simplices at the same time. This becomes a problem if a series of transition *requires* modifying an edge with a strong impact on the simplex count in *multiple* cliques and there is no alternative series of transitions where this situation could be avoided or defused first. In the highest dimensions, where in practical applications double edges are rather rare, this should not happen. In the lower dimensions this is quite likely however, thus we approximate this network effect by assuming that the lower dimensional simplices are embedded a fully-double-edge highest-dimensional *D*-clique. With this assumption, there are $\binom{D-2}{d-2}$ cliques of size $d$ which all share the same (undirected) edge, a multiplier applied to our previous boundary relaxations. The dimension-wise minimum of $\Delta^{(1)}$ and $\Delta^{(2)}$ adjusted by this clustering combinatorial factor should then be wide enough to connect $\mathcal{G}_{s^-}^{s^+}$. We add it the upper end of target boundaries:

$$s_d^{--} = s^- \text{ and } s_d^{++} = \max(s^- + \Delta_d, s^+) \text{ where } \Delta_d = \binom{D-2}{d-2} \cdot \min\left(\Delta_d^{(1)}, \Delta_d^{(2)}\right).$$

## 4 Implementation and Mixing Times

In this section we adopt the perspective of random walks instead of distributions. We start with the rough sketch of the algorithm used for sampling and then discuss the experimental setup used to establish reasonable mixing time.

---
**Algorithm 1:** Restricted Random Walk with Connecting Boundaries
---
    **Input:** Simple directed Graph $G$, desired boundaries $s^-$ and $s^+$,
        number of samples $n$.
    **Optional:** connecting boundaries $s^{--}$ and $s^{++}$, sampling distance $k$.
    **Output:** List of samples from $\mathcal{G}^{s^+}_{s^-}(G)$

[1]  $s \leftarrow$ `count_simplices`$(G)$
[2]  if not supplied, calculate $s^{--}, s^{++}$ and $k$           // see 4.1.5
[3]  `samples` $\leftarrow$ new empty set
    // step 1: MCMC-sampling
[4]  **for** $i \leftarrow 0 \ldots n$ **do**
[5]     $t \leftarrow$ `random_transition`$(G)$               // see 4.1.2
[6]     $G' \leftarrow t(G)$                          // see 4.1.3
[7]     $s' \leftarrow s + s_\Delta$                     // see 4.1.3
[8]     **if** $s^{--} \leq s' \leq s^{++}$ **then**
[9]         $G \leftarrow G'$
[10]       $s \leftarrow s'$
[11]    **if** $i \ \% \ k = 0$ **then**
[12]       `samples.add(G)`
    // step 2: rejection sampling
[13] **for** $G \in$ `samples` **do**
[14]    **if not** $s^- \leq s(G) \leq s^+$ **then**
[15]       `samples.remove(G)`

---

Here the main loop of lines [3]-[12] describes a random walk implementing MCMC-sampling on $\mathcal{G}^{s^{++}}_{s^{--}}(G)$, lines [8]-[10] in particular describe the idea of a restricted walk with resampling, as in case of $G' \notin \mathcal{G}^{s^{++}}_{s^{--}}(G)$ the old graph $G$ is kept, but the counter $i$ is still increased. Lines [13]-[15] describe the subsequent filtering process to only keep the desired graphs in $\mathcal{G}^{s^+}_{s^-}(G)$.

If the connecting boundaries and the sampling distance is correctly set, this algorithm will yield samples from $\mathcal{G}^{s^+}_{s^-}(G)$ whose distribution is arbitrarily close to uniform (see Section 2.2).

## 4.1 Implementation

The implementation is split into two Rust crates: `flag-complex`[4] containing the graph representation and the simplex enumeration code ported from `flagser` [10], and `flag-complex-mcmc`[5], containing code specific to the sampler presented here.

### 4.1.1 Input and Output

Graphs are read in the `.flag` format of `flagser` (ignoring the edge weights) and are written as edge lists in one HDF5 file. An alternative invocation resumes from an automatically saved checkpoint. Invocation examples and Python scripts interpreting the output can be found in the repository[6] containing the setup for the experiments. The sampler is fully deterministic once command

---
[4] https://github.com/TheJonny/flag-complex
[5] https://github.com/TheJonny/flag-complex-mcmc
[6] https://github.com/flomlo/MCMC_directed_flag_complex_sampling_experiments

line supplied integer serving as a seed for the pseudo random number generator is set.

### 4.1.2 Calculating a Transition

At every step, we first choose the kind of transition to apply. The distribution can be selected by a command line option and defaults to 10% SEF, 10% DEM, 60% CP, and 20% CS, values that seemed reasonable but where never optimized. A separate function for each kind of transition returns a list of directed edges to be deleted or inserted.

- Single Edge Flip: A SEF is determined by a uniformly sampled single edge $(i, j)$. We return $(i, j)$ marked for deletion and $(j, i)$ marked for insertion. Here we can safely assume that not all edges are double edges, as there would be only one such graph.

- Double Edge Move: The directed single edge and the undirected double edge are sampled uniformly together with a bit indicating the direction to remove from the double edge. When there is no double edge, we resort to a transition that does nothing.

- Clique Permute: We only consider maximal cliques of the underlying undirected graph. Instead of sampling them uniformly, we first decide the size, weighted by the fifth root of the number of maximal cliques with this size in the graph to increase the chance of larger cliques. Then, for the selected size, a maximal clique is drawn uniformly. The permutation is then chosen uniformly.

- Clique Swap: Two equally sized cliques (with vertex sets $A$ and $B$) are chosen as above, and the bijection required by Definition 3.7 is constructed from three uniformly chosen permutations: One for the intersection of $A$ and $B$, the second for the bijection from $A \setminus B$ to $B \setminus A$, and the third for the bijection from $B \setminus A$ to $A \setminus B$. It is possible that the same cliques are chosen, in which case this forms a Clique Permute.

### 4.1.3 Applying the Transition

With the list of changes generated we then try the transition and check if the simplex count is still inside the given bounds. If not, the transition is reversed. Mutating the graph in-place avoids copying it, as the graph is typically much larger than the transition changeset.

Instead of counting (and therefore enumerating, a computationally costly process) all simplices after each transition we calculate the effect on the simplex count in the neighbourhood `nbhd` of the affected edges (see Definition 2.8). We then calculate the change of the simplex count by restricting the graph to `nbhd`. It follows from Lemma 2.9 that $s(G') = s(G) + s(G'|_{\text{nbhd}}) - s(G|_{\text{nbhd}})$.

Exploiting this locality property is key to making the transition step fast and thus MCMC-sampling computationally feasible.

#### 4.1.4 Initialization and Data Structures

To make the operations stated above efficient we precompute and cache properties only dependent on the underlying undirected graph, as it is invariant with respect to the transitions: we search for the maximal undirected cliques via the Bron-Kerbosch algorithm [21], and save them as arrays of vertex indices.

We also build a hash map from the undirected edges (sorted pairs of vertex indices) to their neighbourhoods saved as an array of vertices according to Definition 2.8. This speeds up calculating the neighbourhood of the transitions.

The current graph is represented both as a dense and a sparse adjacency matrix. The sparse adjacency matrix is an indexable hash set of edges (integer pairs). It provides $\mathcal{O}(1)$ edge insertion, deletion, and uniform sampling. Just like `flagser`, we maintain a dense matrix ($|V|^2$ bits) as it is faster for read only lookups in practice. To uniformly sample double edges quickly, another indexable hash set is maintained. This combined graph representation is implemented in the `EdgeMapGraph` structure in the crate `flag-complex`.

#### 4.1.5 Connecting Boundaries and Sampling Distance

The values for $s^{--}, s^{++}$ are calculated as described in Section 3.3 but can easily be hardcoded as well. The sampling distance is command-line and defaults to $2s_1 \log_2 s_1$, where $s_1$ describes the number of directed edges. This is a reasonable lower bound on sampling distances (compare for mixing times for random walks on hypercubes, see [8]). We dedicate the next section to validating these sampling distances in computational experiments.

### 4.2 Sampling Distance

A rigorous theoretical, yet practically useful result on the mixing time is out of scope for this paper, as it would require extensive theoretical insight into the geometry of the space of directed graphs and simplex counts. Still, the question of what sampling distance to choose is crucial in practice. In this section we describe and conduct two experiments to estimate reasonable defaults used in our implementation and for the experiments described in Section 5: analyzing the distance between two graphs (here the Hamming distance on the adjacency matrix) over time and a test to check for independence of samples.

#### 4.2.1 Hamming Distance over Time

**Definition 4.1** (Hamming distance). Let $G_1 = (V, E_1)$, $G_2 = (V, E_2)$ be two graphs over the same vertex set. The *Hamming distance* is defined by the number of edges unique to one of the graphs

$$d_H(G_1, G_2) = |E_1 \setminus E_2 \cup E_2 \setminus E_1|.$$

This definition is equal to the Hamming distance of the bit strings given by the graphs adjacency matrices.

We measure the development of the distance with respect to the number of Markov Chain transitions $t \in \{0, \cdots, T\}$, where $T \in \mathbb{N}$ is the maximal temporal distance considered. To avoid measuring only the distance to the original starting graph, we sample one very long chain ($893T$) and investigate not only from the starting graph, but also from the offsets $t_0 \in \{0, T, 2T, \cdots\}$. The behaviour of the distance over time is then described by

$$t \mapsto d_H(G_{t_0}, G_{t_0+t}).$$

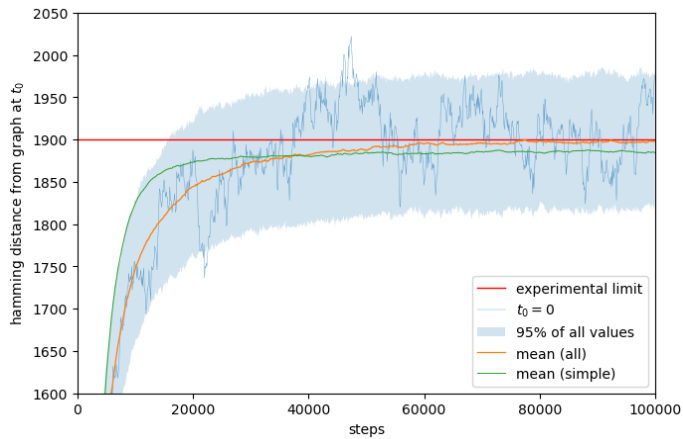

Figure 9: Hamming distance by time, starting at C.Elegans for $t_0 = 0$.

For Figure 9 we started the sampler with the connectome of C. Elegans, target bounds of $\pm 1\%$ and default relaxed bounds and seed. We plot the range of observed values (blue, shaded) and mean (orange line) over all offsets $t_0 \in \{0, \cdots, 892\}$. The behaviour with $t_0 = 0$ (starting graph is C. Elegans) is plotted explicitly (blue line). The red horizontal line is the (experimental) limit behaviour of the Hamming distance. A quick visual analysis seems to suggest that the mean is sufficiently close to the limit after roughly $50000 \simeq 2s_1 \log s_1$ steps.

We further tested the sampler restricted to simple moves as defined in Definition 3.2 (see green line in Figure 9). Initially simple moves move away from the initial state much faster, as they exhibit roughly half the rejection rate compared to complex moves and thus resample less often. But ultimately the green line can be observed to struggle to converge to the limit. We interpret that as a justification of the introduction of the complex moves, even though they take on average 8 times longer to compute per step.

We ran a similar analysis on the BBP layer 1-3 connectome with similar results (not shown).

### 4.2.2 $\chi^2$-Test for Independence

Assume the Markov Chain used by the sampler yields the sequence of graphs $G_t$ with $t \in \mathbb{N}_0$. We want to find the number of steps $k$ required between two graphs $G_t$ and $G_{t+k}$ such that we cannot reject the hypothesis

"$G_t$ and $G_{t+k}$ are stochastically independent".

We check for multiple sampling distances $k$ which are then fixed for each hypothesis test, using a binarizing functions $f : \mathbf{DiGraph} \to \{0, 1\}$ to derive bit sequences $B_n = f(G_{nk})$ with $n \in \mathbb{N}_0$. We now count the *pairs* $(B_{2n}, B_{2n+1}) \in$

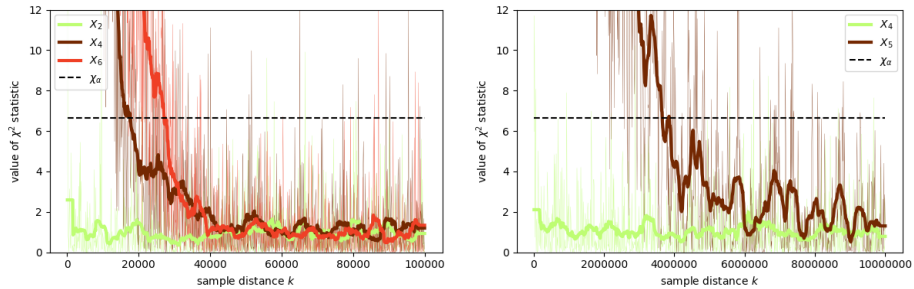

Figure 10: $\chi^2$ Test statistic $X_{d,k}$ using the bit sequences derived from simplex counts in different dimensions. **Thin, colored**: Values from individual tests. **Thick, same color**: smoothed. **Black, dashed**: Acceptance limit $\chi^2_{1-\alpha}$ for $\alpha = 1\%$. If the samples are independent, the probability of $X_k < \chi^2_{1-\alpha}$ is 99%. The dimensions shown here are exemplary, including the longest rejecting ones.
**Left:** Starting with C. Elegans. **Right**: Starting with BBP layers 1-3.

$\{00, 01, 10, 11\}$ of the bit sequences and use Pearson's $\chi^2$ test for independence [3]: if they are independent, the frequency of the pairs should match the frequencies expected by counting single bits. Thus the final test hypothesis is

"$B_{2n}$ and $B_{2n+1}$ are stochastically independent for all applicable $n$".

As binary binning functions we split $s_d(G_t)$ by the median for each dimension $d$.

The test results are shown in Figure 10. We show the values of the test statistic while varying the sampling distance $k$ and the simplex dimension $d$ used for binning. For C. Elegans the data consists of 704964 graphs sampled with distance 100, and we must reject the independence hypothesis for sampling distances smaller than 40000. For BBP layers 1-3 the data consists of 108846 graphs sampled with distance 10000, and we must reject the independence hypothesis for sampling distances smaller than $6 \cdot 10^6$.

This test confirms, both for C. Elegans and BBP layer 1-3, the choice of sampling distances $k \simeq 2s_1 \log_2(s_1)$ presented in Section 4.1.5.

# 5 Application

## 5.1 Homology Analysis for C. Elegans

C. Elegans is a well studied specimen in neuroscience. Its neural wiring is determined completely by its genome and thus it has been studied extensively, leading to an exhaustive map of all its neurons and synapses. We explore the data in the version described in [20], with 279 neurons and 2194 directed synapses.

For our connectomic analysis we ignore any additional information associated with the neurons and synapses (like type, size or location in $3d$-space) and abstract everything away which cannot be represented in a directed graph in the following way: A neuron forms a vertex and a directed chemical synapse with $i$ as its presynaptic and $j$ as its postsynaptic neurons becomes an edge $(i, j)$ in a graph henceforth referred to as $G$.

Using flagser [6, 11] we calculate the simplex count and Betti numbers, with same results as in Figure 1, left. As already discussed in that study, the simplex count and Betti numbers are very distinct to the ones found in comparable Erdős–Rényi graphs.

We pursue a target-relaxation of $\pm 1\%$ and use boundaries as described in Section 3.3.2, with the addition of $s_8^{++} = 10$ as $G$ has cliques of size 9 but no simplices of dimension 8. Thus the simplex counts and boundaries are:

| $d$ | 0 | 1 | 2 | 3 | 4 | 5 | 6 | 7 | 8 |
|---|---|---|---|---|---|---|---|---|---|
| $s^{--}$ | 279 | 2194 | 4276 | 4852 | 4404 | 2681 | 891 | 153 | |
| $s^{-}$ | 279 | 2194 | 4276 | 4852 | 4404 | 2681 | 891 | 153 | |
| $s(G)$ | 279 | 2194 | 4320 | 4902 | 4449 | 2709 | 901 | 155 | |
| $s^{+}$ | 279 | 2194 | 4363 | 4951 | 4493 | 2736 | 910 | 157 | |
| $s^{++}$ | 279 | 2194 | $\infty$ | 4972 | 5124 | 5081 | 2691 | 213 | 10 |

$T$ respects the boundaries $s^{--}, s^{++}$ by resampling instead (around 40% of transitions were accepted) and chooses its moves SEF, DEM, CP, CS with the probabilities [0.1, 0.1, 0.6, 0.2]. For the sampling distance $48704 \simeq 2s_1 log_2(s_1)$ was chosen as motivated in Section 4.2.

Investing 8 days of computation time using a 32 core AMD EPYC 7543 CPU we sampled 10000 samples of $\mathcal{G}_{s^{--}}^{s^{++}}$ for each of the 128 starting seeds. After rejection everything not in $\mathcal{G}_{s^{-}}^{s^{+}}$ 2408 samples remained. This low yield is due to the comparatively small number of 7-simplices in C.Elegans, where $\pm 1\%$ leaves a very small margin for error.

On the remaining 2408 graphs we analyzed the homology using flagser. The result is pretty clear: $G_{\text{C.Elegans}}$ exhibits significantly lower Betti numbers in low dimensions (1 to 3), but higher Betti numbers in high dimensions (6 and 7) than expected on the average of the null model. These differences are significant, often several standard deviations:

| $d$ | 0 | 1 | 2 | 3 | 4 | 5 | 6 | 7 |
|---|---|---|---|---|---|---|---|---|
| mean($b_d$) | 1 | 217 | 283 | 201 | 134 | 92 | 9.8 | 1.3 |
| std($b_d$) | 0 | 6.6 | 15 | 22 | 27 | 21 | 6.6 | 1.4 |
| original | 1 | 183 | 249 | 134 | 105 | 63 | 19 | 5 |
| deviation | | 5.2 | 2.3 | 3.1 | 1.1 | 1.4 | -1.4 | -2.7 |

> Note to the Reviewer: We will make the sampled, filtered graphs available on the long-term scientific data hosting service of the TU Graz. We still need to figure out some details, a footnote with the URL will be added in time.

## 5.2 Homology Analysis for Layers 1-4 of a Rodent Somatosensory Cortex

Computationally, C. Elegans is nowadays seen more as a toy example. To demonstrate that the MCMC-approach is feasible for significantly bigger connectomes as well, we perform the same analysis on the connectome of layers 1-4 of the Blue Brain Project statistical reconstruction of the somatosensory cortex of a rat. With 12519 neurons and 1431930 directed synapses it about half (vertices) to a quarter (edges) of the full reconstruction.

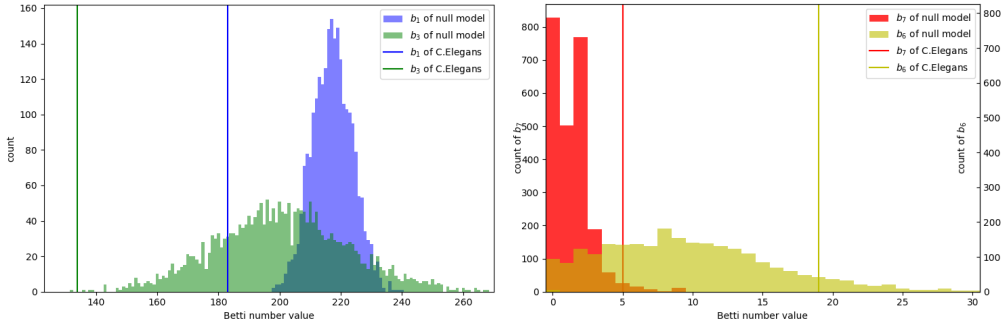

Figure 11: Histogram on the value of the Betti numbers in the null model in dimension 1 and 3 (left) and 6 and 7 (right). Vertical lines show respective values of C.Elegans.

| $d$ | 0 | 1 | 2 | 3 | 4 | 5 | 6 |
|---|---|---|---|---|---|---|---|
| $s^{--}$ | 12519 | 1431930 | 6640109 | 2788378 | 162518 | 1636 | 1 |
| $s^{-}$ | 12519 | 1431930 | 6640109 | 2788378 | 162518 | 1636 | 1 |
| $s(G)$ | 12519 | 1431930 | 6707181 | 2816544 | 164160 | 1653 | 2 |
| $s^{+}$ | 12519 | 1431930 | 6774252 | 2844709 | 165801 | 1669 | 3 |
| $s^{++}$ | 12519 | 1431930 | $\infty$ | 2844709 | 165801 | 2116 | 3 |

We invested the computational power of 32 AMD EPYC 7543 cores over roughly 24 days to sample around 300 samples per core, each sample differs by 114554400 transitions, twice the sampling distance recommended in Section 4.2 to make up for the low transition acceptance rate of around 25%. After rejection-filtering these 9600 samples, 407 remained. We present here the Betti numbers of the original graph as well as mean and standard deviation of the graphs of the null model:

| $d$ | 0 | 1 | 2 | 3 | 4 | 5 |
|---|---|---|---|---|---|---|
| mean($b_d$) | 3 | 16300 | 2672650 | 119456 | 460 | 2 |
| std($b_d$) | 0 | 176 | 12611 | 13715 | 132 | 2 |
| original | 3 | 13748 | 2689400 | 41990 | 70 | 0 |
| deviation | 0 | 14.5 | -1.33 | 5.65 | 2.95 | 1 |

Interestingly, the BBP graph exhibits *significantly less* holes in every dimension except $d = 2$. We can therefore conclude that the geometric spaces described by the flag complex are remarkably flat for graphs with a comparable number of simplices. We hypothesize that this is due to the strong overall flow of information: The majority of synapses go from lower levels to the higher ones.

Note to the Reviewer: Links to samples missing, see above.

# 6 Summary and Outlook

**Summary** The preceding sections began with establishing the need for more refined null models than ER-graphs when investigating directed flag complexes.

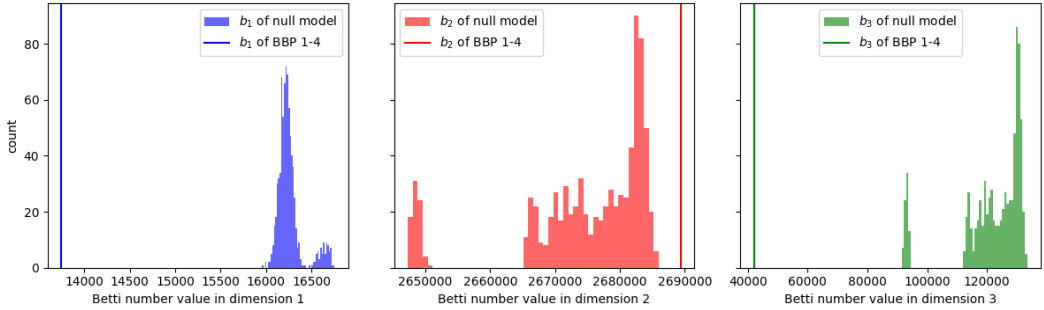

Figure 12: Histogram on the value of the Betti numbers in the null model in dimensions 1 (left), 2 (middle), and 3 (right). Vertical lines show respective values of BBP, layers 1-4.

We decided to first investigate a very specialized null model — it consists of directed graphs with a specified range of simplices in the flag complex and keeps the number of edges and the underlying undirected graph. We then designed, implemented and tested an MCMC-based sampler and made it publicly available. In the special case of oriented graphs (i.e. directed graphs without double edges) we state an elementary graph-theoretical conjecture which, if/once proven, implies correctness and efficiency of our sampler. An investigation of the statistical properties of the homology (i.e. Betti numbers) of two connectomic graphs well-known in Neuroscience then followed and we discovered that they, with some significance, indeed behave extraordinary.

**Implications** We hope that this work leads to a more accurate and quicker assessment of the worth of TDA methods — not only Betti-numbers, but $q$-connectivity and simplicial completions probability — and by that the application of TDA to connectomes and thus Neuroscience in general. In particular, now that we *know* that a variety of connectomes exhibits a topology significantly nonrandom, we should be encouraged to invest more resources into investigating the *how* and *why* of the emergence of holes (or the suspicious lack thereof). We furthermore hope that the insights, approaches and code are of help in developing other models in the area of directed and undirected flag complexes as well as directed (non-flag) simplicial complexes.

**Limitations** From a more applied perspective interested in the graph analysis the null model presented here might be too narrow. By restricting the underlying undirected graph and the number of edges we only allow freedom in the direction of edges and, in the case of non-oriented graphs, the placement of double edges. This corresponds more to concepts of "flow" or "direction". In contrast, the actual "connectivity" is not changed and thus no statements concerning it may be derived from this model.

Simultaneously, the model presented here is too loose for a more rigorous analysis: unlike more traditional approaches like the graph/simplicial configuration model we do not maintain any notion of *directed* vertex or edge participation (the undirected one is trivially maintained due to the fixed undirected backbone). Thus our model does not respect simplex clustering (or a lack thereof) as well as the existing model for undirected simplicial complexes or the related graph configuration model. It could very well be that a null model taking vertex participation into account would not yield comparable results to the one investigated here.

From a computational standpoint there is room for improvement as well: not only very small graphs, where $s^-, s^+$ becomes really narrow, and large graphs, due to the estimated time of $c \cdot s_1 \log(s_1)$ steps, but also (locally) dense graphs where single transitions become very expensive due to the large neighbourhood require substantial computational effort to sample. Even with state-of-the-art CPUs a useful number of samples requires weeks to months of computation time. The most promising and reachable lever is decreasing the rejection rate, i.e. finding better $s^{++}, s^{--}$-bounds.

Both the quest for better boundaries while maintaining irreducibility as well as restricting the model further to respect vertex participation require a substantial better insight into the geometry of $\mathcal{G}_{s^{--}}^{s^{++}}$ in the double-edge case. Incidentally, this is also the biggest theoretical weakness right now: so far we do not have *rigorous* mathematical proofs of irreducibility of the transition matrices. Note that this does not rob this approach of all value — even if the transition matrices used were not irreducible, the statement that the connectomes investigated have unusual topology remains true, although only w.r.t. the connected component. Nonetheless, it would be satisfying *and useful* to be able to formalize, and thus understand, the double-edge case (see Section 3.3.2) as well or even better as the single-edge-case is already understood.

**Potential Research Directions**   As discussed in Section 3.3.1 and Appendix A, the question of connectivity of $\mathcal{G}_{s^{--}}^{\infty}$ for oriented graphs can be stated in a innocuous-looking graph-theoretical way. Proving or refuting Conjecture A.4 is thus an obvious open problem. More generally, the geometry of that space is probably rather tame and worth investigating. Thorough insights would probably enable formal descriptions of mixing time and rejection rate.

The double-edge-case seems considerably more challenging. Simultaneously it is more rewarding, as all practical applications found so for fall into that category. Formalizing, sharpening and proving the hand-wavy assumptions and rules given in Section 3.3.2 is out of scope for the current paper, but should be tackled as soon as possible.

In practical applications though, most of the rejections are due to violations in the highest dimension. From that perspective refining the rules of thumb would be useful enough already. In particular, utilizing knowledge about the undirected graph (i.e. the clique complex, which is fixed!) might be enough to increase the samplers yield by several orders of magnitude.

Once the yield is increased and the necessary computational power is acquired the analysis performed in Section 5 can be applied to bigger and denser connectomes. In particular, the full BBP connectome and the larvae of Drosophila [22] should be addressed.

Similarly we would like to test if *q*-connectedness and the simplex completion probability [16] are as able as Betti-number-analysis in picking up the particularities of actual connectomes w.r.t. the null model here.

Furthermore we hope to combine this work with an (upcoming) one on sampling undirected flag complexes while maintaining clique count (and potentially vertex participation) to construct a sampler which samples directed flag complexes and does not respect the underlying undirected graph. This would form a null model less narrow and artificially restricted.

## Acknowledgements

We would like to express our gratitude for their time, ideas and general help (alphabetic order): Mickaël Buchet, Oliver Friedl, Barbara Giunti, Michael Kerber, Robert Legenstein, Robert Peharz, Henri Riihimäki, Alexander Steinicke, Jean-Gabriel Young. Furthermore we would like to thank two anonymous reviewers for their insightful comments.

## Competing Interest Statement

On behalf of all authors, the corresponding author states that there is no conflict of interest.

## Contribution

Idea and initial concept: FU, developing moves: FU and JK, implementation: JK and FU, mixing time analysis: JK, application: FU, writing: FU and JK.

# A  Converting Graphs into Directed Acyclic Graphs via Single Edge Flips

**Definition A.1** (Directed Graphs)**.** A graph $G = (V, E)$ is

- a simple directed graph iff it does not have 1-cycles, i.e.

$$E \subseteq \{(i,j) \in V \times V |\ i \neq j\}$$

- a simple oriented graph iff it additionally has no 2-cycles (or double edges), i.e.
$$(i,j) \in E \Rightarrow (j,i) \notin E.$$

- a simple oriented complete graph iff additionally for every pair of vertices there exists an edge in one direction or another:

$$(i,j) \in V \times V \Rightarrow \text{ either } (i,j) \in E \text{ or } (j,i) \in E \text{ or } i = j.$$

**Definition A.2** (Single Edge Flip)**.** Let $G = (V, E)$ be a simple oriented graph. Then with $(i,j) \in E$ we define the Single Edge Flip

$$\text{SEF}_{i,j} : G \mapsto G' \text{ where } G' = (V, E') \text{ with } E' = E \setminus \{(i,j)\} \uplus \{(j,i)\}.$$

Note that $G'$ is also a simple oriented graph. If $G$ is complete, so is $G'$.

**Theorem A.3.** *A finite simple oriented complete graph $G = (V, E)$ of $N$ vertices can be transformed into a directed acyclic graph by at most $\frac{|E|}{2} = \frac{N^2 - N}{4}$ SEF such that the number of 3-cycles strictly decreases with each application of SEF.*

*Proof.* The proof works by two loops: The outer loop is an induction over the number of vertices $N$, as with $N = 1$ the statement is trivially true. An inner loop reduces the indegree of a vertex $b$: starting with a (potentially non-unique) vertex $b$ of minimal indegree, we flip incoming edges such that they become outgoing edges to further reduce $\mathrm{indeg}(b)$ until it becomes 0. The vertex $b$ is then a source and the subproblem of $G \setminus \{b\}$ is considered, forming the inductions step $N \mapsto N-1$ of the outer loop. The step where $b$ is removed constructs a total order on $V$: the elements removed earlier (which only have outgoing connections to the remaining graph) are considered to be smaller. This results in a DAG.

It remains to verify that this indegree-reducing SEF is legal, i.e. it strictly reduces the number of 3-cycles in the graph. Let $b \in V$ with $\mathrm{indeg}(b)$ minimal, i.e. $\mathrm{indeg}(b) \leq \mathrm{indeg}(v)$ for all $v \in V$. By assumption $\mathrm{indeg}(b) > 0$, otherwise $b$ would already be a source and we could consider the strictly smaller subproblem of $G \setminus \{b\}$. Let $I = \mathrm{In}(b)$ be the set of all incoming vertices and $O = \mathrm{Out}(b)$ the set of all outgoing vertices of $b$. Now consider $x \in I$: Flipping the edge $(x, b)$ creates as many 3-cycles as $\mathrm{outdeg}_I(x)$, the outdegree of $x$ with respect to $I$, and destroys as many 3-cycles as $\mathrm{indeg}_O(x)$, the indegree of $x$ with respect to $O$. Thus a modification of $G$ via a flip of $(x, b)$ strictly decreases the number of 3-cycles and thus adheres to our restriction iff $\mathrm{indeg}_O(x) > \mathrm{outdeg}_I(x)$. We observe (justifications below):

$$\begin{aligned}
\mathrm{indeg}_O(x) &= \mathrm{indeg}(x) - \mathrm{indeg}_I(x) \\
&= \mathrm{indeg}(x) - (|I| - 1 - \mathrm{outdeg}_I(x)) \\
&\geq \mathrm{indeg}(b) - \mathrm{indeg}(b) + 1 + \mathrm{outdeg}_I(x) \\
&= \mathrm{outdeg}_I(x) + 1
\end{aligned}$$

As $G$ is complete, $V = \{b\} \cup I \cup O$ and thus $\mathrm{indeg}(x) = \mathrm{indeg}_O(x) + \mathrm{indeg}_I(x)$, leading to the first line. Another implication of completeness is the relation of indegree and outdegree: as, between to *different* vertices, there must be exactly one edge in one direction or another, $\mathrm{indeg}(x) = |V| - 1 - \mathrm{outdeg}(x)$. This statement still holds true when restricted to $I$, thus the identity $\mathrm{indeg}_I(x) = |I| - 1 - \mathrm{outdeg}(x)$ leads to the second line. The identity $|I| = \mathrm{indeg}(b)$ leads, together with the initial assumption that $b$ has minimal indegree, to the lower bound proposed in the third line. This proves conformity of the inner loop of our requirement to reduce the number of 3-cycles each step.

Note that once an edge $(x, b)$ has been flipped, it will not be flipped again: $b$ will be the unique element of minimal degree in the next step, and $x$ lies now in $\mathrm{Out}(b)$. Once $b$ has become a source, all connections with $b$ will no longer be considered. Thus each edge must be flipped at most once. Furthermore, an element with minimal indegree may have at most $\frac{N}{2}$ incoming edges, leading to the additional factor of $\frac{1}{2}$ in the number of steps required. □

Does this Theorem also hold in the case of non-complete graphs? No.

If one weakens the requirement of a *strictly* monotonic decrease in the number of 3-cycles to just a monotonic decrease in the number of 3-cycles (i.e each SEF may not increase the number of 3-cycles) then the result still seems to be true: The graph can be transformed into a DAG while flipping each edge at most once.

We tested this hypothesis on ER-graphs configured with ($N = 30$, $p = 0.8$), ($N = 60$, $p = 0.4$) and ($N = 300$, $p = 0.1$) on 150 million graphs each and did not find any counterexample. However, despite ongoing effort, the following conjecture still remains unproven:

**Conjecture A.4.** *A finite simple oriented graph $G = (V, E)$ of $N$ vertices can be transformed into a directed acyclic graph by at most $|E|$ applications of SEF such that the number of 3-cycles does not increase with each individual step.*

# References

[1] Yael Artzy-Randrup and Lewi Stone. "Generating uniformly distributed random networks." In: *Phys Rev E Stat Nonlin Soft Matter Phys* 72.5 Pt 2 (Nov. 2005), p. 056708. DOI: 10.1103/PhysRevE.72.056708.

[2] Yazan N Billeh et al. "Systematic Integration of Structural and Functional Data into Multi-scale Models of Mouse Primary Visual Cortex." In: *Neuron* 106.3 (May 2020), 388–403.e18. DOI: 10.1016/j.neuron.2020.01.040.

[3] The SciPy community. *scipy.stats.chi2_contingency*. https://docs.scipy.org/doc/scipy-1.11.2/reference/generated/scipy.stats.chi2_contingency.html. [accessed 04-September-2023]. 2023.

[4] Steven J Cook et al. "Whole-animal connectomes of both Caenorhabditis elegans sexes". In: *Nature* 571.7763 (2019), pp. 63–71.

[5] Bailey K. Fosdick et al. "Configuring random graph models with fixed degree sequences". In: *SIAM Review* 60.2 (2018), pp. 315–355. ISSN: 0036-1445. DOI: 10.1137/16M1087175.

[6] giotto-tda. *pyflagser is a python API for the flagser C++ library*. https://github.com/giotto-ai/pyflagser. 2019.

[7] Chad Giusti et al. "Clique topology reveals intrinsic geometric structure in neural correlations". In: *Proceedings of the National Academy of Sciences of the United States of America* 112.44 (2015), pp. 13455–13460. ISSN: 0027-8424. DOI: 10.1073/pnas.1506407112.

[8] David A Levin and Yuval Peres. *Markov chains and mixing times*. Vol. 107. American Mathematical Soc., 2017.

[9] Sahil Loomba et al. "Connectomic comparison of mouse and human cortex". In: *Science* 377.6602 (2022), eabo0924. DOI: 10.1126/science.abo0924. eprint: https://www.science.org/doi/pdf/10.1126/science.abo0924. URL: https://www.science.org/doi/abs/10.1126/science.abo0924.

[10] Daniel Lütgehetmann. *flagser*. https://github.com/luetge/flagser. 2017-2021.

[11] Daniel Lütgehetmann et al. "Computing persistent homology of directed flag complexes". In: *Algorithms (Basel)* 13.1 (2020), Paper No. 19, 18. DOI: 10.3390/a13010019.

[12] Henry Markram et al. "Reconstruction and simulation of neocortical microcircuitry". In: *Cell* 163.2 (2015), pp. 456–492.

[13] Alessandro Motta et al. "Dense connectomic reconstruction in layer 4 of the somatosensory cortex." In: *Science* 366.6469 (Nov. 2019). DOI: 10.1126/science.aay3134.

[14] Theodore Papamarkou et al. "A random persistence diagram generator". In: *Statistics and Computing* 32.5 (2022), Paper No. 88, 15. ISSN: 0960-3174. DOI: 10.1007/s11222-022-10141-y.

[15] Michael W. Reimann et al. "Cliques of Neurons Bound into Cavities Provide a Missing Link between Structure and Function." In: *Frontiers Comput. Neurosci.* 11 (2017), p. 48. URL: http://dblp.uni-trier.de/db/journals/ficn/ficn11.html#ReimannNSTPCDLH17.

[16] Henri Riihimäki. *Simplicial q-connectivity of directed graphs with applications to network analysis*. Feb. 2022. arXiv: 2202.07307v1 [math.AT]. URL: http://arxiv.org/abs/2202.07307v1;%20http://arxiv.org/pdf/2202.07307v1.

[17] Louis K Scheffer et al. "A connectome and analysis of the adult Drosophila central brain". In: *Elife* 9 (2020), e57443.

[18] Brad Theilman, Krista Perks, and Timothy Q. Gentner. "Spike Train Coactivity Encodes Learned Natural Stimulus Invariances in Songbird Auditory Cortex". In: *Journal of Neuroscience* 41.1 (2021), pp. 73–88. ISSN: 0270-6474. DOI: 10.1523/JNEUROSCI.0248-20.2020. eprint: https://www.jneurosci.org/content/41/1/73.full.pdf. URL: https://www.jneurosci.org/content/41/1/73.

[19] Florian Unger, Jonathan Krebs, and Michael Müller. "Simplex closing probabilities in directed graphs". In: *Computational Geometry* 109 (Feb. 2023). DOI: 10.1016/j.comgeo.2022.101941.

[20] Lav R Varshney et al. "Structural properties of the Caenorhabditis elegans neuronal network". In: *PLoS Comput Biol* 7.2 (2011), e1001066.

[21] Wikipedia. *Bron–Kerbosch algorithm — Wikipedia, The Free Encyclopedia*. http://en.wikipedia.org/w/index.php?title=Bron%E2%80%93Kerbosch%20algorithm&oldid=1110744344. [accessed 27-September-2022]. 2022.

[22] Michael Winding et al. "The connectome of an insect brain". In: *Science* 379.6636 (2023), eadd9330.

[23] Jean-Gabriel Young et al. "Construction of and efficient sampling from the simplicial configuration model". In: *Physical Review E* 96.3 (Sept. 2017). ISSN: 2470-0053. DOI: 10.1103/physreve.96.032312. URL: http://dx.doi.org/10.1103/PhysRevE.96.032312.